\documentclass[12pt]{article}
\usepackage{amsfonts}
\newtheorem{theo}{Theorem}[section]
\newtheorem{pro}[theo]{Proposition}

\newcommand{\ra}{\rightarrow}
\newcommand{\varep}{\varepsilon}

\title{ Nash 
type inequalities for fractional powers of non-negative self-adjoint 
operators}
\author{A. Bendikov\\
{\small Department of Mathematics, Cornell University}\\
{\small Ithaca, NY, 14853-4201, USA}
\and
P. Maheux
\footnote
{Research partially supported by {\sl European commission} (TMR 
1998-2001 Network {\sl Harmonic Analysis}).}\\
{\small D\'epartement de Math\'ematiques.U.M.R 6628.
MAPMO}\\
{\small Universit\'e d'Orl\'eans B.P 6759}\\
{\small 45 067 Orl\'eans Cedex 2, France}}

\begin{document}
    
\maketitle

\begin{abstract}
Assuming that a Nash type inequality is satisfied by a non-negative 
self-adjoint operator
$A$,  we prove a Nash type inequality for the fractional 
powers $A^{\alpha}$ of $A$. Under some
assumptions, we give  ultracontractivity bounds for the
  semigroup $(T_{t,{\alpha}})$ generated by $-A^{\alpha}$.
\end{abstract}


\noindent

\section{Introduction} \label{secintro}
 \setcounter{equation}{0}

Let $(T_{t})$ be a symmetric submarkovian semigroup acting on 
$L^2(X,{\mu})$ with ${\mu}$ a ${\sigma}$-finite measure on $X$ and 
let $(-A,{\cal D})$ be its generator. The following theorem is 
known and due to Varopoulos and Carlen, Kusuoka and Stroock (see 
\cite{VSC} Thm II.5.2 and references therein).
\begin{theo}
    For  $n>2$, the following conditions are equivalent :
\begin{equation}\label{soba}
  \mid\mid f \mid\mid_{2n/n-2}^2\leq C(Af,f),\quad \forall f\in {\cal 
D}; 
\end{equation}

\begin{equation}\label{nasha}
  \mid\mid f \mid\mid_{2}^{2+4/n}\leq C_{1}(Af,f)\mid\mid f 
  \mid\mid_{1}^{4/n},\quad \forall f\in {\cal D}\cap L^1(X,{\mu}); 
\end{equation}

\begin{equation}\label{rega}
  \mid\mid T_{t} \mid\mid_{1\rightarrow \infty } \leq C_{2}t^{-n/2} 
  ,\quad \forall t>0 . 
 \end{equation}
\end{theo}
In particular, using subordination, $(\ref{nasha})$ 
implies that for all ${\alpha}\in (0,1)$ :

\begin{equation}\label{nashapha}
  \mid\mid f \mid\mid_{2}^{2+4{\alpha}/n}\leq C_{3}(A^{\alpha}f,f)\mid\mid 
f 
  \mid\mid_{1}^{4{\alpha}/n},\quad \forall f\in {\cal 
D}(A^{\alpha})\cap L^1(X,{\mu}).
\end{equation}

In  \cite{C} and \cite{D}, an equivalence of the type  
 (\ref{nasha})-(\ref{rega}) was proved in
greater generality under some assumptions on the function
$t\rightarrow \mid\mid T_{t} \mid\mid_{1\rightarrow \infty } $. In particular 
in \cite{C} the following condition $(D)$ is used.
\\

\noindent
{\bf Definition.} {\em A differentiable function $m$:
${\mathbb R}_{+}\longrightarrow {\mathbb R}_{+}$ satisfies 
condition $(D)$ if the function
$M(t):=-{\log}\, m(t)$ is such that :
$$
\forall t>0,\quad \forall u\in [t,2t],\qquad M^{\prime}(u)\geq c 
M^{\prime}(t)$$
for some constant $c>0$.}

Let $m_{1},m_{2}$ be two functions from $]0,+\infty[$ to itself; we 
shall say that  $m_{1} \preceq m_{2}$   if there 
exist $C_{1},C_{2}>0$ 
such that $m_{1}(t)\leq C_{1} m_{2}(C_{2} t)$, and  that $m_{1}$, $m_{2}$  
are equivalent $(m_{1}\simeq m_{2})$ if $m_{1}\preceq m_{2}$ and $m_{2}\preceq
m_{1}$. Note that if $m_{i},\;i=1,2,$ are decreasing differentiable bijections 
satisfying $(D)$ and if ${\Theta}_{i}(x)=-m_{i}^{\prime}(m_{i}^{-1}(x))$, then 
$m_{1}\simeq m_{2}$ if and only if ${\Theta}_{1}\simeq {\Theta}_{2}$.
In the two following statements, the inequalities will be written 
modulo equivalence of functions.

\begin{theo}
Let $m$ be a decreasing $C^1$ bijection of ${\mathbb R}_{+}$ satisfying $(D)$
and set ${\Theta}(x)=-m^{\prime}(m^{-1}(x))$. Then the following 
conditions are equivalent :

\begin{equation}\label{tetanash}
{\Theta}( \mid\mid f \mid\mid _{2}^2)\leq (Af,f),\quad \forall f\in 
{\cal D}(A),\;\;  \mid\mid f \mid\mid_{1}=1.
\end{equation}
\begin{equation}
 \mid\mid T_{t} \mid\mid_{1\rightarrow \infty } \leq m(t)
   ,\quad \forall t>0 . 
\end{equation}
 \end{theo}

We consider the following question : Assume that 
$A$ satisfies the Nash type inequality (\ref{tetanash}).
What kind of Nash inequality is satisfied by the operator $A^{\alpha}$, 
the fractional power of $A$?
In what follows, it will be convenient to write  (\ref{tetanash})
in the equivalent form (see \cite{BM})
 
 \begin{equation}\label{caln}
  \mid\mid f \mid\mid _{2}^2 B\left(   \mid\mid f \mid\mid 
_{2}^2\right)\leq (Af,f),\quad \forall f\in 
{\cal D}(A),\;\;  \mid\mid f \mid\mid_{1}=1,
\end{equation}
where
 \begin{equation}\label{be}
     B(x)=\sup_{t>0}(t\log x+tM(1/t)),\quad  M(t)=-\log m(t).
  \end{equation}
In particular, $x\rightarrow B(x)$ is a non-decreasing function
satisfying the following property
$$\lim_{x \rightarrow\infty}\frac{B(x)}{\log x}=+\infty.$$ 
In some cases (non-ultracontractive semigroups) an inequality similar to 
(\ref{caln}) 
can be proved with $x\rightarrow B(x)$ being 
non-decreasing but not necessary obtained from a function $m$ by (\ref{be}). 
For instance, the function $x\rightarrow B(x)$ 
with 
$B(x)=\log x$ may be relevant (see  Section  \ref{orn}) .  We state our result 
with a very weak  assumption on $x\rightarrow B(x)$ in order to take into 
account such cases.
The main result of this note is the following theorem 
 
 \begin{theo}\label{main}
 Let $(X,{\mu})$ be a measure space with  ${\sigma}$-finite measure 
${\mu}$.
 Let $A$ be a non-negative self-adjoint operator with  
domain ${\cal D}(A) \subset L^2(X,{\mu})$. 
Suppose that the semigroup   $T_{t}=e^{-tA}$ acts as a contraction on 
$L^1(X,{\mu})$ and satisfies the following Nash type inequality 

\begin{equation}\label{Nashone}
\mid\mid f \mid\mid_{2}^2 B\left(\mid\mid f \mid\mid_{2}^2\right)
 \leq (Af,f),\quad
\forall f \in
{\cal D}(A),\;\;\mid\mid f \mid\mid_{1} = 1,
\end{equation} 
where $B:[0,+\infty[\ra   [0,+\infty[ $ is
a non-decreasing function which tends to infinity at infinity. 
Then, for any ${\alpha}>0$,  the following 
Nash type inequality holds
 \begin{equation}\label{Nalphanash}
\mid\mid f \mid\mid _{2}^2  \left[ B\left(   
\mid\mid f \mid\mid 
  _{2}^2\right)\right]^{\alpha}  \leq (A^{\alpha}f,f),\quad \forall f\in 
{\cal D}(A^{\alpha}), \;\; \mid\mid f \mid\mid _{1}=1.
\end{equation}
\end{theo}
\noindent
{\bf Remark.} Thus, if the function $x\ra B(x)$ corresponds 
by (\ref{caln}) to the 
operator  $A$ then the function $x\ra \left[ B(x)\right]^{\alpha}$ 
 corresponds to the operator $A^{\alpha}$.
 The function $x\ra x^{\alpha}, x\geq 0$, $0\leq {\alpha}\leq 1$ is a 
particular case of 
 so-called Bernstein function (see \cite{BF}). The importance of 
Bernstein 
 functions comes from the following property.
 If $-A$ is a Markov generator then for any Bernstein function $g$, 
  $-g(A)$ is again a Markov generator. More precisely, $-g(A)$ 
  generates Markov semigroup $(T_{t}^g)$ given by the following formula 
  $$
  T_{t}^g=\int_{0}^{+\infty} T_{s}\;d{\mu}_{t}^g(s), \quad t>0,
  $$
  where $(T_{s})$ is the Markov semigroup generated by $-A$ and 
  $({\mu}_{t}^g)_{t\geq 0}$ is the one-sided-stable convolution 
  semigroup on ${\mathbb R}_{+}$ (the subordinator) defined  uniquely by 
  its Laplace transform
  $$
 \int_{0}^{+\infty}e^{-xs}\;d{\mu}_{t}^g(s)=e^{-tg(x)},\quad x>0.
 $$
 In view of Theorem \ref{main}  one may wonder if the Nash 
 inequality (\ref{Nashone}) for $-A$ implies Nash inequality for 
 $-g(A)$ in the form
 
   \begin{equation}\label{Noconj}
\mid\mid f \mid\mid_{2}^2 g\circ B\left(\mid\mid f \mid\mid_{2}^2 \right) \leq
\left(g(A)f,f\right),\;\;\forall f\in
  {\cal D}(g(A)),\;\;  \mid\mid f\mid\mid_{1}=1.
   \end{equation}
\noindent   
In general, the answer is not known. For instance, although we strongly 
suspect that for  
a minimal Bernstein function $g:x\ra 
1-e^{-ax}, a>0$, (\ref{Nashone}) does not implies (\ref{Noconj}), we have 
no proof of this fact at the present writing. 
It would be interesting to describe the set of Bernstein functions 
for which we can pass from Nash inequality (\ref{Nashone})
 to Nash inequality (\ref{Noconj}). Theorem \ref{main} states that this 
set contains all power 
functions $x\ra x^{\alpha},\;0<{\alpha}\leq 1$.

\section{Proof of Theorem \ref{main} and related 
results}\label{proofsect}\setcounter{equation}{0}

In this section, we will prove  Theorem \ref{main} 
in three steps. We prove (\ref{Nalphanash}) with ${\alpha}=1/2$.
Then we iterate the result of step 1 to prove 
(\ref{Nalphanash}) for all ${\alpha}_{n}$ of the form 
${\alpha}_{n}=1/2^n$, $n\in {\mathbb N}$.
We give a convexity argument which will  allow  us to  
conclude for 
$0<{\alpha}<1$ and also for $ {\alpha}\geq1$.

Before  embarking on  the proof, we need some preparations.
 Let $A$ be a non-negative self-adjoint operator on $L^2(X,{\mu})$,
${\mu}$ is a ${\sigma}$-finite measure. Since $A$ is non-negative, its 
spectral decomposition has the form
$$A=\int_{0}^{+\infty}{\lambda}\,dE_{\lambda}.$$
In particular, the semigroup $T_{t}=\int_{0}^{\infty} 
e^{-t{\lambda}}\,dE_{\lambda}$  generated by $-A$ satisfies 
$\mid\mid T_{t} \mid\mid_{2\ra 2}\leq 1$ for all $t>0$.
The fractional power $A^{\alpha}$ of $A$ is defined  by 
the formula
$$(A^{\alpha}f,f)=\int_{0}^{+\infty}{\lambda}^{\alpha}\,d(E_{\lambda}f,f)$$
  on   the domain
${\cal D}(A^{\alpha})=\{ f\in L^2  :
\int_{0}^{+\infty}{\lambda}^{2{\alpha}}\,d(E_{\lambda}f,f) < +\infty\}$.
The operators $A^{\alpha}$ are   non-negative self-adjoint operators. 
The contraction  semigroup generated by $-A^{\alpha},0<{\alpha}<1,$ can 
be expressed 
in the form
$$T_{t,{\alpha}}=\int_{0}^{+\infty} T_{s}\;d{\mu}_{t}^{\alpha}(s),$$
where $({\mu}_{t}^{\alpha})$ is the one-sided  ${\alpha}$-stable semigroup 
on 
${\mathbb  R}_{+}$. This semigroup 
can be characterized by its Laplace transform 
$$\int_0^{+\infty} e^{-t{\lambda}}\; d{\mu}_{t}^{\alpha}(s) 
=e^{-t{\lambda}^{\alpha}},\quad {\lambda}>0. $$
We also denote  $T_{t,1/2}$ by $P_t$ and call  
$(P_t)$ the Poisson semigroup associated to $A$. 

\subsection{The case ${\alpha}=1/2$ }
\begin{theo}\label{square}
Let $A$ be a non-negative self-adjoint operator such that $T_t=e^{-tA}$ 
acts as a contraction on $L^1$ for all $t> 0$.
Assume that there exists $B: {\mathbb R}^+\ra {\mathbb R}^+$, 
non-decreasing and 
such that
\begin{equation}\label{nashb}
\mid\mid f\mid\mid_2^2B(\mid\mid f\mid\mid_2^2)\leq (Af,f),\quad \forall 
f\in {\cal D}(A),\;\;\mid\mid f\mid\mid_{1}=1.
\end{equation}
Then, for all ${\varepsilon}\in (0,1)$,  
\begin{equation}\label{nashbdemi}
(1-{\varep}^2)^{1/2}
\mid\mid f\mid\mid_2^2\left[ B({\varep}\mid\mid 
f\mid\mid_2^2)\right]^{1\over2}
\leq (A^{1/2}f,f),\quad\forall f\in {\cal 
D}(A^{1/2}),\;\; \mid\mid f\mid\mid_{1}=1.
\end{equation}
\end{theo}
\noindent
{\bf Proof.} Let $g\in {\cal D}(A^{1/2})$ and $\mid\mid g\mid\mid_{1}\leq 
1$. Set $f=P_t g$. 
Then $f\in  {\cal D}(A^{n})$, for all $n\geq 1$. 
 The semigroups $(P_t)$ and $(T_t)$ are related by the subordination 
formula

$$P_t g=\int_0^{+\infty} {\mu}_{t}^{1/2}(s) T_s g\, ds=
\frac{1}{\sqrt{\pi}}\int_0^{+\infty}\frac{e^{-u}}{\sqrt{u}} T_{t^2/4u} 
f\,du.
 $$
It follows that $(P_t)$ is a contraction 
semigroup on $L^1$ and in particular
$$\mid\mid f\mid\mid_{1}=\mid\mid P_t g\mid\mid_{1}\leq \mid\mid 
g\mid\mid_{1}\leq 1.$$
Since $f\in {\cal D}(A)$, the relation $Af=AP_t g=\frac{d^2}{dt^2} 
P_t g$ 
holds true. We  apply (\ref{nashb}) with $f=P_t g$  
\begin{equation}\label{eqpt}
\mid\mid P_t g\mid\mid_2^2B(\mid\mid P_t g\mid\mid_2^2)\leq (AP_t g, P_t g).
\end{equation}
Set ${\phi}(t)=\mid\mid P_t g\mid\mid_2^2$, then

$$\frac{d}{dt}{\phi}(t)=\dot{\phi}(t)=-2(A^{1/2}P_t g,P_t g)$$
and
$$ \frac{d^2}{dt^2}{\phi}(t)=\ddot{\phi}(t)=4(AP_t g,P_t g).$$
The inequality (\ref{eqpt}) can be written in the form
\begin{equation}\label{eqphi}
4{\phi}(t)B({\phi}(t))\leq {\ddot{\phi}}(t),\;\; t>0.
\end{equation}
 Multiplying both sides in
(\ref{eqphi}) by $-\dot{\phi}\geq 0$, we obtain 
$$
-4[{\phi}^2(t)]^{\prime}B({\phi}(t))\leq 
-[{\dot{\phi}}^2]^{\prime}(t),\quad \forall t>0.
$$
Fix $T>0$ and integrate this inequality over $[0,T]$ to obtain
$$
-4\int_0^T B({\phi}(s))[{\phi}^2(s)]^{\prime}\,ds\leq 
-\int_0^T [{\dot{\phi}}^2]^{\prime}(s)\,ds.
$$
The right hand side is clearly bounded by 
$[{\dot{\phi}}]^2(0)=4(A^{1/2}g,g)^2$
for all $T>0$.
To deal with the left hand side, set $v(s)={\phi}^2(s)$. Then it takes 
the form
$$-4\int_0^TB(\sqrt{v(s)})  v^{\prime}(s) \, ds=4
\int_{v(T)}^{v(0)} B(\sqrt{x})\,dx.$$
Thus finally we get the following inequality
\begin{equation}\label{eqint}
\int_{v(T)}^{v(0)} B(\sqrt{x})\,dx\leq (A^{1/2}g,g)^2.
\end{equation}
Let us assume that 
\begin{equation}\label{spe}
\lim_{T\ra +\infty} \mid\mid P_T \; g\mid\mid_2=0.
\end{equation}
Below, we will see  how to reduce the general case to this one. 
In the inequality (\ref{eqint}), we take the limit as ${T\ra +\infty}$ 
and obtain 
$$
\int_{0}^{v(0)} B(\sqrt{x})\,dx\leq (A^{1/2}g,g)^2.
$$
Let ${\varep}\in (0,1)$. Since $B$ is non-decreasing,
$$
(1-{\varep}^2)v(0)
B({\varep}\sqrt{v(0)})\leq 
\int_{{\varep}^2v(0)}^{v(0)} B(\sqrt{x})\,dx\leq (A^{1/2}g,g)^2
$$
and finally,
\begin{equation}\label{eqdemi}
(1-{\varep}^2)^{1/2}
\mid\mid g\mid\mid_2^2\left[ B\left({\varep}\mid\mid 
g\mid\mid_2^2\right)\right]^{1\over 2}
\leq (A^{1/2}g,g).
\end{equation}
This proves the theorem under the assumption (\ref{spe}).    
To consider the general case, define the operator
$A_{\rho}=A+{\rho}I$, ${\rho}>0$. $A_{\rho}$ is non-negative and 
self-adjoint. 
It  also satisfies  (\ref{nashb}). The property 
$$
\lim_{T\ra +\infty} \mid\mid e^{-T\sqrt{A+{\rho}I}}
\mid\mid_2=0
$$
follows by  spectral theory. 
We apply  the inequality (\ref{eqdemi})   with $A_{\rho}$ instead of $A$.
 Since the left hand side of (\ref{eqdemi}) is independant of ${\rho}>0$, 
 we can pass to
the limit as ${\rho}\ra 0$ in 
(\ref{eqdemi}). 
The proof is now complete.

\subsection{Iteration}

\begin{pro}
Under the same assumptions as  in Theorem \ref{square}, for all $n\in 
{\mathbb N}^*$ there
exists ${\alpha}_n,{\beta}_n>0$ such that,
\begin{equation}
{\alpha}_n \mid\mid f\mid\mid_2^2\left[ 
B({\beta}_n\mid\mid f\mid\mid_2^2)\right]^{1/{2^n}}
\leq ( A^{1/{2^n}}f,f),\quad
\forall f\in {\cal D}(A^{1/{2^n}}),\;\; \mid\mid f\mid\mid_{1}=1 
\end{equation}
\end{pro}
{\bf Proof.} We apply Theorem \ref{square} and induction on $n$.

\subsection{The convexity argument}
We have already proved  (\ref{Nalphanash}) for 
 ${\alpha}={\alpha}_n=1/{2^n}$,
$n\in {\mathbb N}^*$. To conclude that (\ref{Nalphanash}) holds true
for all ${\alpha}\in (0,1)$ we need the following auxiliary result.

\begin{pro}\label{convarg}
    Let ${\Lambda}:{\mathbb R}^{+}\longrightarrow {\mathbb R}^{+}$ be 
a non-decreasing function. 
    Assume that $A$ is a non-negative self-adjoint operator that  
satisfies the inequality 
$$
\mid\mid f\mid\mid_2^2{\Lambda}(\mid\mid f\mid\mid_2^2)\leq 
(Af,f),\quad\forall 
f\in {\cal D}(A), \mid\mid f\mid\mid_{1}=1. 
$$
Then for any convex non-decreasing
function ${\Phi}\geq 0$
$$
\mid\mid f\mid\mid_2^2 {\Phi}\,\circ\,{\Lambda}(\mid\mid f\mid\mid_2^2)\leq 
({\Phi}(A )f,f),\quad \forall 
f\in  {\cal D}({\Phi} (A )), \mid\mid f\mid\mid_{1}=1. 
$$
\end{pro}
{\bf Proof.}
By renormalisation, $f\ra f/ \mid\mid f \mid\mid_{1}$, we have
$$
\mid\mid f\mid\mid_2^2{\Lambda}(\mid\mid f\mid\mid_2^2/\mid\mid 
f\mid\mid_1^2)
\leq 
\int_0^{+\infty} {\lambda}\, d(E_{\lambda}f,f).
$$
For a fixed $f$ denote $d{\nu}(\lambda)= d(E_{\lambda}f,f)$.
Assume that $\mid\mid f\mid\mid_2=1$, then $\nu$ is a probability measure. 
Since ${\Phi}$ is   convex non-decreasing function,  Jensen's inequality 
yields
$$
{\Phi}\,\circ\,{\Lambda}(1/\mid\mid f\mid\mid_1^2)
\leq 
\int_0^{+\infty} {\Phi}({\lambda} )\, d(E_{\lambda}f,f),
$$
that is
$$
{\Phi}\,\circ\,{\Lambda}(1/\mid\mid f\mid\mid_1^2)
\leq 
({\Phi} (A )f,f),\quad\mid\mid f\mid\mid_{2}=1.
$$
This obviously gives the result.

\subsection{End of the proof}

Let $0<{\alpha}<1 $ be fixed and choose $n\in {\mathbb N}^*$ such that
${\alpha}_n=1/{2^{n}}\leq {\alpha}$. We   have 
$$
a_{{\alpha}_n}\mid\mid f\mid\mid_2^2 \left[ B(
b_{{\alpha}_n}\mid\mid f\mid\mid_2^2)\right]^{{\alpha}_n}
\leq (A^{{\alpha}_n}f,f),\quad
\forall f\in {\cal D}(A^{{\alpha}_n}),\;\;
\mid\mid f\mid\mid_{1}=1. 
$$
Choose 
${\Phi}(t)=t^{{\alpha}/{\alpha}_n}$ and let
${\Lambda}(t)=a_{{\alpha}_n}\left[ B(b_{{\alpha}_n}t)\right]^{{\alpha}_n}$.
 Since ${\alpha}/{\alpha}_n \geq 1$, ${\Phi}$ is a non-decreasing 
convex function.
Moreover,
$${\Phi}\,\circ\,{\Lambda}(t)=a_{{\alpha}}\left[ B(b_{{\alpha}}t) 
\right]^{\alpha}, $$
where
$a_{\alpha}=(a_{{\alpha}_n})^{{\alpha}/{{\alpha}_n}}$
, $b_{{\alpha}}=b_{{\alpha}_n}$.
For $f\in {\cal D}(A^{\alpha})$, ${\Phi}(A^{{\alpha}_n})f=A^{\alpha}f$
and Proposition \ref{convarg}  yields the result
$$
a_{\alpha} \mid\mid	f \mid\mid _{2}^2  
\left[ B(b_{\alpha}\mid\mid	f\mid \mid _{2}^2)\right]^{\alpha}
\leq (A^{\alpha}f,f),\quad\forall f\in {\cal D}(A^{{\alpha}}),\;\;
\mid\mid f\mid\mid_{1}=1.
$$
This finishes the proof of Theorem \ref{main} for $0<{\alpha}\leq 1 $.
In the case ${\alpha}>1$, we just apply Proposition \ref{convarg}.

\subsection{Some generalizations }\label{secnashreg} 
We want to enlarge the class of functions treated in Theorem  \ref{main}.
Recall the notion of  regularly varying function 
(see \cite{BGT}).
Function  ${\Phi}$ defined on $[0,+\infty[$ is said to be 
regularly varying function of index ${\alpha}$ if 
for any ${\lambda}\geq 1$,
$$\lim_{x\ra 
\infty}\frac{{\Phi}({\lambda}x)}{{\Phi}(x)}={\lambda}^{\alpha}.
$$
In the case ${\alpha}=0$, ${\Phi}$ is called a slowly varying function.
Any regularly varying function of index ${\alpha}$ can be represented 
in the form  ${\Phi}(x)=x^{\alpha}{\ell}(x)$, where ${\ell}$ is a slowly 
varying function.
\\

\noindent
{\bf Examples.} The following
functions illustrate the definition of  regular  variation of index $\alpha$:
$ x \rightarrow cx^{\alpha},  cx^{\alpha}(\log x)^{\beta},
cx^{\alpha}(\log x)^{\beta}(\log \log x)^{\gamma},
cx^{\alpha}\exp\left[ ( \log x)^{\delta}\right],
$
where 
$-\infty <{\alpha}, {\beta}, {\gamma}  <+\infty$ and $0<{\delta}<1$.    

\begin{theo}\label{regu}
Let $A$ be a non-negative self-adjoint operator and ${\Phi} :{\mathbb 
R}^+ \longrightarrow {\mathbb R}^+$
a regularly varying function of index ${\alpha}>0$. Assume that   
there exists  $B:{\mathbb R^+} \longrightarrow {\mathbb R}^+$ such 
that $B(x)\nearrow\infty$ as $x\nearrow\infty$ and
\begin{equation}\label{nashn}
\mid\mid f\mid\mid_2^2 B(   \mid\mid f\mid\mid_2)
\leq (Af,f),\quad \forall f\in {\cal D}(A),\;\; \mid\mid f\mid\mid_1=1.
\end{equation}
Then there exist $c,a>0$ such that 
 \begin{equation}\label{nashnphi}
c\mid\mid f\mid\mid_2^2 {\Phi}\,\circ\,B(   \mid\mid f\mid\mid_2)
\leq ({\Phi}(A)f,f),\quad \forall f\in {\cal D}({\Phi}(A)),\;\;  \mid\mid 
f\mid\mid_1=1,
\mid\mid f\mid\mid_2\;\geq a.
\end{equation}
\end{theo}
The proof of Theorem \ref{regu} is based on the following auxiliary
results.

\begin{pro}\label{phi2}
Let ${\Phi} :{\mathbb R}^+ \longrightarrow {\mathbb R}^+$ be such that for 
some $N\geq 1$, the
function $\varphi : x \longrightarrow {\Phi}(x^N)$ is eventually  increasing and 
convex.
Then (\ref {nashn})
implies (\ref {nashnphi}).
\end{pro}
To prove proposition \ref{phi2} we apply Theorem \ref{main},  the 
convexity argument of Proposition 
 \ref{convarg} and the following relation  
$$ {\Phi}(A)f={\varphi}(A^{1\over N})f\;,\quad \forall  f\in 
{\cal D}(A)\cap {\cal D}({\Phi}(A)).$$
 
\begin{pro}\label{var}
Let ${\varphi} :{\mathbb R}^+ \longrightarrow {\mathbb R}^+$ be a 
regularly varying function of  index
${\alpha}>1$. Then there exists ${\Phi} :{\mathbb R}^+ \longrightarrow 
{\mathbb R}^+$
which is  regularly varying of  index
${\alpha}$, increasing, convex and such that 
$$\lim_{x\longrightarrow +\infty} \frac{{\Phi}(x)}{{\varphi}(x)}=1.$$
\end{pro}
{\bf Proof.} 
The function ${\varphi}$ can be represented in the form 
$${\varphi}(x)= x^{\alpha}{\ell}(x),Ê\quad x>0,$$
where ${\ell}$ is slowly varying and non-negative.
Define the following functions 
$${\tilde{\varphi}}(x)={\alpha}({\alpha}-1) x^{{\alpha}-2}{\ell}(x),$$
$$ {\Phi}(x)=\int_0^x d{\tau}\int_0^{\tau} {\tilde{\varphi}}(s)\;ds.$$
It follows that $  {\Phi}^{\prime}$ and  $  {\Phi}^{{\prime}{\prime}}$
are non-negative. Hence ${\Phi}$ is increasing and convex. The function 
${\tilde{\varphi}}$ is regularly varying of index ${\alpha}-2$.
By Feller's theorem,  
$$\int_0^{\tau}{\tilde{\varphi}}(s)\;ds 
\sim
 \frac{1}{{\alpha}-1}{\tau}
{\tilde{\varphi}}({\tau})={\alpha}  {\tau}^{{\alpha}-1}{\ell}({\tau}),
\quad {\tau}\rightarrow +\infty$$
and
$$
\int_0^x\left(\int_0^{\tau}{\tilde{\varphi}}(s)\;ds \right)\;d\tau
\sim
 \frac{1}{{\alpha}}x\left({\alpha}x^{{\alpha}-1}{\ell}(x)\right)={\varphi}(x),
 \quad x\rightarrow +\infty.
$$
This finishes the proof.
\\

\noindent
{\bf Proof of Theorem \ref{regu}}
For any fixed $N>\frac{1}{{\alpha}}$,
the function $x\rightarrow {\Phi}(x^N)$ is regularly varying of index 
${\alpha}^{\prime}={\alpha}N>1$. By Lemma \ref{var}  there exists 
${\tilde{\Phi}}$ regularly varying of index ${\alpha}^{\prime}$, 
increasing,
convex and such that ${\Phi}(x^N)\sim {\tilde{\Phi}}(x), 
\quad x\rightarrow +\infty$.
Set $H(x):= {\tilde{\Phi}}
(x^{1\over N})$. Then there exists $a_{1}>0$ such that
$$ 
{1\over 2} H(x)
\leq
{\Phi}(x)
\leq 
2H(x)
,\quad \forall  x\geq a_{1}.
$$
We apply now Proposition  \ref{phi2}. For any $f\in L^2$ with 
$\mid\mid f\mid\mid_1=1$
we have 
$$
({\Phi}(A)f,f)=\int_{0}^{\infty}{\Phi}({\lambda})\,d(E_{{\lambda}}f,f)
= 
\int_{a_{1}}^{\infty}{\Phi}({\lambda})\,d(E_{{\lambda}}f,f)
+
\int_{0}^{a_{1}}{\Phi}({\lambda})\,d(E_{{\lambda}}f,f)
$$
$$
\geq {1\over 2}\int_{a_{1}}^{\infty}H( {\lambda} )\,d(E_{{\lambda}}f,f)
-c_{1}\mid\mid f\mid\mid_{2}^2
\;\;
\geq
{1\over 2}\int_{0}^{\infty}H( {\lambda} )\,d(E_{{\lambda}}f,f)
-c_{2}\mid\mid f\mid\mid_{2}^2
$$
$$
={1\over 2}(H(A)f,f)-c_{2}\mid\mid f\mid\mid_{2}^2
\;\;\geq 
{1\over 2}c_{0}\mid\mid f\mid\mid_{2}^2\,H{\circ}B(\mid\mid f\mid\mid_{2})
-c_{2}\mid\mid f\mid\mid_{2}^2.
$$
Since $H$ and $B$ approach to infinity as $x\ra \infty$ we can find 
$a>a_{1}$ such that for $x\geq a, B(x)\geq a_{1}$ and 
${\Phi}\circ B(x)\geq 4c_{2}/c$. Hence for  $f\in L^{1}\cap L^{2}$, such 
that 
$\mid\mid f\mid\mid_{1}=1, \mid\mid f\mid\mid_{2}\geq a$
$$
({\Phi}(A)f,f)\geq {1\over 4}c_{0}\,{\Phi}{\circ}B(\mid\mid 
f\mid\mid_{2})-c_{1}\mid\mid f\mid\mid_{2}^2\;\;\geq
{1\over 8} c_{0}\,{\Phi}{\circ}B(\mid\mid f\mid\mid_{2}).
$$
This finishes the proof of the theorem.

\section{Contraction properties of the semigroup
$T_{t,{\alpha}}$}\setcounter{equation}{0}

Let  $(T_t)$ be  a semigroup acting on all $L^p, 1<p<\infty$.  
$(T_t)$ is said to be ultracontractive if for every $t>0$, the operator
$T_t$ can be extended to a bounded operator from  $L^1$  to $L^{\infty}$. 
That is,
there exists a non-decreasing function $m$ from 
${\mathbb R}_{+} $
to itself such that 
$$ \mid\mid T_t \mid\mid_{1\ra \infty}\leq m(t),\qquad t>0.$$

\noindent
$(T_t)$ is said to be hypercontractive if there exists $t>0$ such that 
$T_t$ is a bounded operator from  $L^2$  to $L^4$. 
See \cite{G}.
 
 In the following theorem, all inequalities will be understood  
 in the sense of equivalent functions. See Section \ref{secintro}. 

\begin{theo}\label{exp}
Let $A$ be a non-negative self-adjoint operator such that the 
semigroup $T_{t}=e^{-tA},\;t>0,$ acts as a contraction semigroup on $L^{1}$.
\begin{enumerate}
\item
The following properties are equivalent
\begin{enumerate}
\item
There exits ${\gamma}>0$, such that for any $t>0$,
\begin{equation}\label{expo}
  \mid\mid T_{t} \mid\mid_{1\rightarrow \infty } \leq 
   e^{t^{-{\gamma}}}
 \end{equation}
\item
The following Nash inequality  holds
\begin{equation}\label{nashexpo}
  \mid\mid f \mid\mid_{2}^{2}\left[ \log_{+}(\mid\mid f 
  \mid\mid_{2}^{2} \right]^{1+1/{\gamma}} \leq 
 ( Af,f),\quad f\in{\cal D}(A),\;  \mid\mid f 
  \mid\mid_{1}=1
 \end{equation}
\end{enumerate}
 \item
 Assume that the equivalent properties  1(a) and 1(b) 
 hold. Let $0<{\alpha}\leq 1$,
 then the following inequality holds 
\begin{equation}\label{nashexpoalpha}
  \mid\mid f \mid\mid_{2}^{2}\left[ \log_{+}(\mid\mid f 
  \mid\mid_{2}^{2} \right]^{{\alpha}(1+1/{\gamma})} \leq 
( A^{\alpha}f,f)
,\quad  f\in{\cal D}(A^{\alpha}),\;  \mid\mid f 
  \mid\mid_{1}=1
 \end{equation}
\noindent
 In particular, let ${\alpha}_{c}=\frac{{\gamma}}{{\gamma}+1}$, then
\begin{enumerate}
    \item
 If $ {\alpha}>{\alpha}_{c}$,
 then $T_{t,{\alpha}}=e^{-A^{\alpha}t}$ is ultracontractive and
$  \mid\mid T_{t,{\alpha}}\mid\mid_{1\rightarrow \infty } \leq 
  e^{t^{-{\beta}}}$,
 where ${\beta}=\frac{ {\alpha}_{c} }{{\alpha} - {\alpha}_{c} }$.
 \item
 If   $ {\alpha}\leq{\alpha}_{c}$, then $(T_{t,{\alpha}})$ may  not be
 ultracontractive. See Section 4.
\item If $ {\alpha}={\alpha}_{c}$, and $-A$ is a Markov generator,  
the following logarithmic 
  Sobolev inequality   holds. There exists $C>0$ such that 
   \begin{equation}\label{gross}
   \int f^{2} \log \left(\frac{f}{ \mid\mid f \mid\mid_{2}} \right)\; 
d{\mu} \leq C\left[\left( A^{\alpha}f,f\right)+ 
   \mid\mid f \mid\mid_{2}^{2}\right],
   \quad f\in{\cal D}(A^{\alpha})
   \end{equation}
 In particular, $ (T_{t,{\alpha}})$ is hypercontractive.
\end{enumerate}
 \end{enumerate}
 \end{theo}
{\bf Proof.}  
Statement 1  is a consequence of Theorem 1.2. Statement
2(a) follows from Statement 1 and Theorem \ref{main};
${\beta}=\frac{ {\alpha}_{c} }{{\alpha} - {\alpha}_{c} }$ is 
the result of the  integration of the Nash inequality (\ref{nashexpoalpha}),
see \cite{D}.
For 2(b) we refere to Theorem \ref{ben} (1) below. 
In order to consider the  case ${\alpha}={\alpha}_{c}$,
we need the following   result  from \cite{BM}, see also \cite{BCLS}.
\begin{pro}\label{tron}
 Suppose that $({\cal E},{\cal D}) $  is a quadratic form in 
 $L^2(X,{\mu})$ which  satisfies the 
 following conditions 
\begin{enumerate}
    \item
 For any non-negative $f\in {\cal D}$,
 $f_{k}=(f-2^k)^+\wedge 2^k\in {\cal D}$ for 
 all $k\in {\mathbb Z}$,
 \item 
 $\,\sum_{k\in {\mathbb Z}} {\cal E}(f_{k})\leq 
 {\cal E}(f)$,
\item
  For any non-negative $f\in {\cal D},\mid\mid f \mid\mid_{1}\leq 1$,
$\, \mid\mid f \mid\mid _{2}^2\log \mid\mid f \mid\mid _{2}
 \leq {\cal E}(f)$.   
\end{enumerate}
Then  
there exists a constant
$C>0$ such that
  \begin{equation}\label{grossnash}
  \int f^{2} \log \left(\frac{f}{ \mid\mid f \mid\mid_{2}}\right)\; d{\mu} 
\leq C\left[{\cal E}(f)+ 
  \mid\mid f \mid\mid_{2}^{2}\right],\quad  f\in {\cal D},\;f\geq 0.
 \end{equation}
\end{pro}

For the sake of completeness, we give the proof of this statement.
Let   $f\in~{\cal D},f\geq~0$.
 Without lost of generality, we assume that 
$\mid\mid f \mid\mid_{2}=1$. Let $f_{k}$ be  as above, then $\mid\mid 
f_{k} \mid\mid _{1}<\infty$. For all
$k\in {\mathbb Z}$ we have 
$$ \mid\mid f_{k} \mid\mid _{2}^2\log (\mid\mid f_{k} \mid\mid _{2}
/\mid\mid f_{k} \mid\mid _{1})
 \leq {\cal E}(f_{k}).$$
Since $ \mid\mid f  \mid\mid_{2}= 1$  we have  
$\mid\mid f_{k} \mid\mid _{2}  
/\mid\mid f_{k} \mid\mid _{1}\geq 2^{k}$. 
Indeed,
$$
\mid\mid f_{k} \mid\mid _{1}=\int_{\{f_k>0\}} f_kd{\mu}
\leq {\mu}(f_k>0)^{1/2} \mid\mid f_{k} \mid\mid _{2}
={\mu}(f>2^k)^{1/2} \mid\mid f_{k} \mid\mid _{2}
$$
$$
\leq 2^{-k} \mid\mid f \mid\mid _{2}
 \mid\mid f_{k} \mid\mid _{2}
\leq
2^{-k} \mid\mid f_{k} \mid\mid _{2}.
$$   
Markov inequality 
and the inequality above imply
$$(2^k)^2{\mu}(f\geq 2^{k+1})\log (2^{k-1})\leq {\cal E}(f_{k}).$$
Let $A_{k}=\{2^{k+1}\leq f<2^{k}\}$, then we have
$$\int_{X}  f^{2} \log  f \; d{\mu}
=\sum_{k\in {\mathbb Z}}\int_{A_{k}}f^{2} \log  f \; d{\mu}
\leq
\sum_{k\in {\mathbb Z}}
(2^{k+2})^2{\mu}(f\geq 2^{k+1})\log (2^{k+2}).
$$
This  yields,
$$\int_{X}  f^{2} \log  f \; d{\mu}
\leq 16\sum_{k\in {\mathbb Z}} {\cal E}(f_{k})+12\log 2\left(
\sum_{k\in {\mathbb Z}}(2^{k+1})^2{\mu}(f\geq 2^{k+1})\right).$$
We conclude by (2) and by the fact that the last sum is comparable to 
$\mid\mid f \mid\mid_{2}^2$.
\\

\noindent
{\bf Proof of Theorem \ref{exp} (c)} Since $-A$ is a Markov generator,
${\cal E}(f):=(Af,f)$ is a Dirichlet form. Thus (1) and (2) of 
Proposition \ref{tron} hold true. Property (3) follows from the Nash 
inequality (\ref{nashexpoalpha}) with 
${\alpha}=\frac{\gamma}{(1+{\gamma})}$. Thus, by Proposition 
\ref{tron}, the logarithmic Sobolev inequality (\ref{grossnash}) holds. This 
implies 
 hypercontractivity of $(T_{t,{\alpha}})$, see [G,\,Theorem 3.7].

\section{Invariant Dirichlet forms on the infinite 
dimensional torus}\setcounter{equation}{0}

In this section, we consider the case where the measure space 
$(X,{\mu})$ is the infinite dimensional torus ${\mathbb T}^{\infty}$,
the product of countable many 
copies of ${\mathbb T}= {\mathbb R}/ 2\pi{\mathbb Z} $. The 
topology on ${\mathbb T}^{\infty}$ is the product topology  generated 
by cylindric sets. We regard ${\mathbb T}^{\infty}$ as a 
compact connected abelian group equipped with its 
(normalized) Haar measure ${\mu}$ and will
focus on invariant strictly local Dirichlet forms $({\cal E},{\cal F})$ 
on the group ${\mathbb T}^{\infty}$. 
All the examples below are taken from \cite{B} and \cite{BSC} and the 
aim of this section is to illustrate the results of Sections 1,2 and 3.
We assume that both $\cal F$ and $\cal E$ are invariant under the action of
translations on functions. Any such Dirichlet form can be described by a 
symmetric 
nonnegative definite matrix $A=(a_{i,j})$ so that the associated Dirichlet
form is given on smooth cylindric functions by the formula
$${\cal E}(f,f)=\int_{{\mathbb T}^{\infty}}\sum_{i,j}a_{i,j}{\partial}_{i} f
{\partial}_{j} f\; d{\mu}.$$
Yet another characterisation of $\cal E$ is that the $L^2$-generator $L$ 
associated
to $\cal E$ on smooth cylindric functions is given by the formula
$$ Lf =-\sum_{i,j}a_{i,j}{\partial}_i{\partial}_{j}f.$$
Because of translation invariance, the associated semigroup $T_t:=T_t^A$ is
given by convolution with a Gaussian semigroup of measures $({\mu}_t^A)$,
that is $T_t^A f={\mu}_t^A*f$. See Heyer's book \cite{He} for background on
convolution semigroups of measures on locally compact groups.

\subsection{The product semigroup $T_t^A$}\label{diagcase}
Assume that $A$ is a diagonal matrix with diagonal entries 
$a_{k,k}:=a_{k}$. In this case, 
$\displaystyle{\mu}_{t}^A=\otimes_{1}^{\infty}{\eta}_{a_{k}t}$ is a 
product-measure, where  $({\eta}_{s})_{s>0}$ is the standard Gaussian 
convolution semigroup on the torus  ${\mathbb T}$. Since the operators $T_t^A$ 
act 
as convolutions one can show that the semigroup $(T_t^A)$ is ultracontractive 
if and only if the measures ${\mu}_t^A$ are absolutely continuous w.r.t. 
${\mu}$ and admit continuous densities $x\ra {\mu}_t^A(x)$. In this case
$$
\mid\mid T_t^A \mid\mid_{L_1\ra L_{\infty}}={\mu}_t^A(e),\quad e=(0,0,\cdots).
$$  
Define the following function 
$$
N_{A}(s)= \sharp\{ k:a_{k}\leq s\}, \quad s>0.
$$
Then, the measures ${\mu}_{t}^A$ are absolutely continuous w.r.t. ${\mu}$  
if and 
only if $\log N_{A}(s)=o(s)$ as $s\nearrow  \infty$. In this case,  
the densities  $x\ra {\mu}_{t}^A(x)$  are continuous functions. 
Moreover, if we assume that 
$N_{A}$ varies regularly of index 
${\gamma}>0$, then there exists $C_{\gamma}>0$ such that
$$
\log\mid\mid T_{t}^A \mid\mid_{L_{1}\ra L_{\infty}}=\log{\mu}_{t}^A(e)
\sim
C_{\gamma}N_{A}\left(\frac{1}{t}\right),\quad t\searrow 0.
$$
In particular, if $N_{A}(s)\sim s^{\gamma}$ as $s\nearrow \infty$ then
$$
\log\mid\mid T_{t}^A \mid\mid_{L_{1}\ra L_{\infty}}
\sim
C_{\gamma}t^{-\gamma},\quad t\searrow 0.
$$

\subsection{Contraction properties of the semigroup $T_{t,{\alpha}}^A$}

We now apply  the results of Section \ref{diagcase}
to the semigroup $T_{t,{\alpha}}^A$ 
generated by the operator $-(L_{A})^{\alpha}$, $0<\alpha <1$. This 
clearly will illustrate the results of Section 3. In what follows, 
we assume that $N_{A}(s)\sim s^{\gamma}$ as $s\ra \infty$. Hence ${\mu}_t^A$ 
is absolutely continuous w.r.t. $\mu$ and admits a continuous density 
${\mu}_t^A(x)$ for all $t>0$. Moreover  
condition (\ref{expo}) holds in very precise form
$$
 \log\mid\mid T_{t}^A \mid\mid_{L_{1}\ra L_{\infty}}
 \sim
 C_{\gamma}t^{-{\gamma}},\quad t\searrow 0.
 $$
Because of the 
subordination relation,
$$
T_{t,{\alpha}}^A f={\mu}_{t,{\alpha}}^A*f,
$$
where  $({\mu}_{t,{\alpha}}^A)_{t>0}$  is a convolution semigroup of 
probability 
measures ${\mu}_{t,{\alpha}}^A$ given by  the formula (see 
Section \ref{proofsect})
$$
{\mu}_{t,{\alpha}}^A=\int_{0}^{\infty} {\mu}_{s}^A\;d{\mu}_{t}^{\alpha}(s).
$$
Since ${\mu}_{s}^A$  is absolutely continuous w.r.t. ${\mu}$ for all 
$s>0$, ${\mu}_{t,{\alpha}}^A $ is absolutely continuous w.r.t. ${\mu}$ for 
all 
$t>0$ and admits a lower 
semi-continuous density $x\ra {\mu}_{t,{\alpha}}^A(x)$. By symmetry
$$
\mid\mid T_{t,{\alpha}}^A \mid\mid_{L_{1}\ra L_{\infty}}=
{\mu}_{t,{\alpha}}^A(e), \quad t>0.
$$
These observations and well-known asymptotic properties of the 
${\alpha}$-subordinator $({\mu}_{t}^{\alpha})_{t>0}$ (see \cite{Z}) give 
the following result. 

\begin{theo}\label{ben}
    Assume that $ N_{A}(s)\sim s^{\gamma}, {\gamma}>0$.
   Let ${\alpha}_{c}=\frac{\gamma}{\gamma+1}$ be the critical 
   exponent (see Theorem \ref{exp} (2)). Then 
 \begin{enumerate}
\item
$\forall {\alpha}\in ]0,{\alpha}_{c}[, \forall t>0$,
  $$
\mid\mid T_{t,{\alpha}}^A \mid\mid_{L_{1}\ra L_{\infty}}
={\mu}_{t,{\alpha}}^A(e)=+\infty. 
$$
 \item 
 $\forall {\alpha}\in ] {\alpha}_{c},1[, \forall t>0$,
  $$
\mid\mid T_{t,{\alpha}}^A \mid\mid_{L_{1}\ra 
L_{\infty}}={\mu}^A_{t,{\alpha}}(e)
< +\infty 
$$
and 
$$
 \log\mid\mid T_{t,{\alpha}}^A \mid\mid_{L_{1}\ra L_{\infty}}
 =\log{\mu}_{t,{\alpha}}^A(e)
 \sim
 C_{\alpha,\gamma}t^{-{\alpha}_{c}/({\alpha} -{\alpha}_{c})},\quad 
t\searrow 0.
 $$
\item 
 ${\alpha}={\alpha}_{c}$. There exists $t=t_{\gamma}>0$ such that
  \begin{enumerate}
 \item
 $\forall t\in ]0,t_{{\gamma}}[$,
 $$
 \mid\mid T_{t,{\alpha}_{c}}^A \mid\mid_{L_{1}\ra L_{\infty}}
 ={\mu}_{t,{\alpha}_{c}}^A(e)=+\infty, 
$$
\item
$\forall t\in ] t_{{\gamma}},\infty[$,
 $$\mid\mid T_{t,{\alpha}_{c}}^A \mid\mid_{L_{1}\ra L_{\infty}}
 ={\mu}_{t,{\alpha}_{c}}^A(e) < +\infty.  
$$
 \end{enumerate}
\end{enumerate}
\end{theo} 

\noindent
{\bf Remark.} Because we assume that $N_{A}(s)\sim s^{\gamma}$, 
$s\ra \infty$, the generator $-L_{A}$ of the semigroup $(T_t^A)$ 
satisfies the following Nash 
inequality
$$
\mid\mid  f\mid \mid_{2}^2\left(\log_{+} \mid\mid  f\mid 
\mid_{2}^2\right)^{1+1/{\gamma}}
\leq 
(L_{A}f,f),\quad \mid\mid  f\mid \mid_{1}\leq 1,
$$
which in this special case becomes sharp, i.e. for some  sequence 
$\{f_{n}\}$ 
such that $\mid\mid  f_{n}\mid \mid_{1}\leq 1$ and 
$\mid\mid  f_{n}\mid \mid_{2}\ra \infty$ as $n\ra \infty$
the LHS and  the RHS of the inequality above are comparable. 
This implies sharpness of the Nash inequality for $ (L_{A})^{\alpha}$

\begin{equation}\label{optim}
\mid\mid  f\mid \mid_{2}^2\left(\log_{+} \mid\mid  f\mid 
\mid_{2}^2\right)^{{\alpha}(1+1/{\gamma})}
\leq 
\left( (L_{A})^{\alpha} f,f \right),\quad \mid\mid  f\mid \mid_{1}\leq 1.
\end{equation}
In particular, if $ {\alpha}<{\alpha}_{c}=\frac{\gamma}{\gamma+1}$,
$(L_{A})^{\alpha}$ does not satisfies the log-Sobolev 
inequality (\ref{grossnash}). 
Indeed, the log-Sobolev inequality for $(L_{A})^{\alpha}$ would imply the Nash 
inequality of the form 
$$
\mid\mid  f\mid \mid_{2}^2
\log \mid\mid  f\mid \mid_{2}^2\leq C
\left( (L_{A})^{\alpha}f,f\right), \quad \mid\mid  f\mid \mid_{1}\leq 1,
$$
for some $C>0$. This is not possible since (\ref{optim}) is sharp.
Hence, by (\cite{G}, Theorem 3.7 ) 
$(T_{t,{\alpha} }^A)_{t>0}$,  ${\alpha}<{\alpha}_{c}\;$, is 
not a hypercontractive semigroup.

\section{Ornstein-Uhlenbeck semigroup} \label{orn}
 \setcounter{equation}{0}

As an example where Theorem \ref{main} applies 
we consider   Ornstein-Uhlenbeck semigroup $(T_{t})$ on  
$L^{2}({\gamma}_{n},{\mathbb R}^n)$, where
 ${\gamma}_n$ is the standard Gaussian measure on ${\mathbb R}^n$

$$ d{\gamma}_n(x)=(2\pi)^{-n/2} e^{-\mid x\mid^2/2}dx\;. $$
\noindent
The Dirichlet form associated to $(T_{t})$ is defined as
$$
{\cal E}(f)=\int_{{\mathbb R}^n} \mid {\nabla}f \mid^2\;d{\gamma}_n, 
$$
where $\mid {\nabla}f \mid^2$ is the square of the 
gradient of $f$. Then the generator $A$ of the form ${\cal E}$ is 
represented on 
${\cal C}_c^{\infty}({\mathbb R}^n)$ by the following equality

$$-A={\Delta}+x{\nabla}.$$
According to (\cite{G}, Example 4.2) $A$ satisfies the following 
log-Sobolev inequality
\begin{equation}\label{lsou}
\int f^2\log \left(f/ \mid\mid f\mid\mid_2\right)\;d{\gamma}_n
\leq
\int \mid {\nabla}f \mid^2\;d{\gamma}_n =(Af,f),\quad \forall f\in{\cal 
C}_c^{\infty}({\mathbb R}^n).
\end{equation}
Thus, $(T_{t})$ is hypercontractive, but not ultracontractive. Indeed, if 
$f(x)=x_{i}$ then $f\in L^{1}$, but $T_{t}f=f\notin L^{\infty}$. 
Inequality (\ref{lsou}) implies

\begin{equation}\label{nou}
\mid\mid f\mid\mid_2^2\log  \mid\mid f\mid\mid_2 
\leq (Af,f),
\quad f\in {\cal D}(A),\;\; \mid\mid f\mid\mid_1= 1.
\end{equation}
\noindent
Since $\mid\mid f\mid\mid_1\leq \mid\mid f\mid\mid_2$, the LHS of 
this inequality is 
non-negative. Theorem \ref{main} implies that  for all 
$0<{\alpha}<\infty $, the following Nash inequality holds

\begin{equation}\label{noualpha}
\mid\mid f\mid\mid_2^2\left( \log  \mid\mid f\mid\mid_2 \right)^{\alpha}
\leq (A^{\alpha}f,f),
\quad f\in {\cal D}(A^{\alpha}),\;\; \mid\mid f\mid\mid_1= 1.
\end{equation}

We claim that for any $0<{\alpha}<1$ the semigroup  $(T_{t,{\alpha}})$
is not hypercontractive. Indeed, hypercontractivity of $(T_{t,{\alpha}})$
would imply the inequality

\begin{equation} \label{false}
\mid\mid f\mid\mid_2^2\left( \log  \mid\mid f\mid\mid_2^{2} \right) 
\leq (A^{\alpha}f,f),
\quad f\in {\cal D}(A^{\alpha}),\;\; \mid\mid f\mid\mid_1= 1.
\end{equation}
\noindent
Then, by Theorem \ref{main}, the inequality (\ref{false}) implies the 
following Nash inequality

\begin{equation} \label{falsealp}
\mid\mid f\mid\mid_2^2\left( \log  \mid\mid f\mid\mid_2^{2} 
\right)^{1\over \alpha}
\leq (A f,f),
\quad f\in {\cal D}(A ),\;\; \mid\mid f\mid\mid_1= 1.
\end{equation}
\noindent
Since ${1\over \alpha}>1$, the Nash inequality (\ref{falsealp}) 
implies ultracontractivity of $(T_{t})$

$$\mid\mid T_{t} \mid\mid_{1\rightarrow \infty}\leq
e^{t^{-{\gamma}}},\;\; {\gamma}=\frac{{\alpha}}{1-{\alpha}}.$$
Contradiction. This proves 
the claim.
\\
The reasons given above imply the following more general result.
 
 \begin{pro}
     Let $-A$ be a symmetric Markov generator. Assume that the 
     semigroup $(T_{t})$ generated by $-A$ is not ultracontractive. 
     Then, for any $0<{\alpha}<1$, the semigroup $(T_{t,{\alpha}})$ 
     generated by $-(A^{\alpha})$ is not hypercontractive.
     \end{pro}
\noindent
{\bf Acknowledgements.}  We would like to thank L. Saloff-Coste from 
Cornell University for fruitfull discussions and valuable comments.

\end{document}